\documentclass{compositio}
\usepackage{amsmath,amsthm,amssymb,amscd}

\theoremstyle{plain}

\newtheorem{cor}{Corollary}[section]

\newtheorem{thm}{Theorem}[section]
\theoremstyle{remark}
\newtheorem{rem}{Remark}[section] 

\theoremstyle{definition}
\newtheorem{defi}{Definition}[section]
\newcommand{\bb}[1]{\mbox{$\mathbb{#1}$}}

\begin{document}
\bibliographystyle{plain}


\title{$ $\\$ $\\$ $\\
A general intersection formula for Lagrangian cycles}

\author{J\"{o}rg Sch\"{u}rmann}
\email{jschuerm@math.uni-muenster.de}
\address{Westf. Wilhelms-Universit\"{a}t\\ 
SFB 478
Geometrische Strukturen in der Mathematik\\
Hittorfstr.~27\\ 48149 M\"{u}nster\\ Germany}



\begin{abstract}
We prove a generalization to the context of real geometry of an intersection formula
for the vanishing cycle functor, which in the complex context is due to Dubson, L\^{e},
Ginsburg and Sabbah (after a conjecture of Deligne). It is also a generalization of similar
results of Kashiwara-Schapira, where these authors work with a suitable assumption
about the micro-support of the corresponding constructible complex of sheaves.
We only use a similar assumption about the support of the corresponding characteristic cycle 
so that our result
can be formulated in the language of constructible functions and Lagrangian cycles.
\end{abstract}

\maketitle


\section*{Introduction}  
In this paper we give a proof (of a generalization to the context of {\em o-minimal
structures} on the real field, or to {\em analytic (Nash) geometric categories}) of 
the following intersection formula:

\begin{thm} \label{thm:int1}
Let $M$ be an m-dimensional real analytic manifold and 
$f: M\to \bb{R}$ a subanalytic
$C^{2}$-function. Consider on $M$ a bounded subanalytically constructible 
complex
${\cal F}$ of sheaves of vector spaces (over a base-field $k$), with finite 
dimensional stalks ${\cal F}_{x}$ ($x\in M$). Suppose that the intersection
of $\sigma_{f}:=\{(x,df_{x})\in T^{*}M | x\in M\}$ and the support 
$|CC({\cal F})|$ of the characteristic cycle of ${\cal F}$ is 
contained in a
compact subanalytic subset $I\subset T^{*}M$, with $K:=\pi(I)\subset \{f=0\}$.
Then one has
\begin{equation} \label{eq:int1}
\chi\bigl(R\Gamma(K,R\Gamma_{\{f\geq 0\}}{\cal F})\bigr) = 
\sharp\bigl(\; [df(M)]\cap [CC({\cal F})]\; \bigr) \:.
\end{equation}
\end{thm}

Here we use the following notations:
\begin{enumerate}
\item $\chi$ is the usual {\em Euler characteristic} (note that 
$R\Gamma(K,R\Gamma_{\{f\geq 0\}}{\cal F})$ has finite dimensional cohomology,
since $K:=\pi(I)$, with $\pi: T^{*}M\to M$ the natural projection, is a
compact subanalytic subset). 
\item $[CC({\cal F})] \in H^{m}_{|CC({\cal F})|}(T^{*}M,\pi^{-1}or_{M})$ is 
the {\em characteristic cycle} of ${\cal F}$ (as in \cite[chap. IX]{KS}, and see the
next section), with $or_{M}$
the {\em orientation sheaf} on $M$. 
\item $[df(M)]\in H^{m}_{\sigma_{f}}(T^{*}M,or_{T^{*}M/M})$ corresponds to
$1\in H^{0}(M,k_{M})\simeq
H^{m}_{\sigma_{f}}(T^{*}M,or_{T^{*}M/M})$ (as in
\cite[def. 9.3.5]{KS}, with $\pi^{!}\bb{Z}_{M}\simeq or_{T^{*}M/M}[m]\;$).
\item The {\em intersection number} $\sharp$ is defined by
\end{enumerate}
\begin{displaymath} \begin{CD}
H^{m}_{\sigma_{f}}(T^{*}M,or_{T^{*}M/M})\times 
H^{m}_{|CC({\cal F})|}(T^{*}M,\pi^{-1}or_{M}) 
@>\cup >> H^{2m}_{c}(T^{*}M,or_{T^{*}M}) 
@> tr >> \bb{Z} \:.
\end{CD} \end{displaymath}


Here we use $or_{T^{*}M/M}\otimes_{Z} \pi^{-1}or_{M} = or_{T^{*}M}$,
and the above cup-product is the composition of the usual cup-product with
support and the natural maps 
\[H^{2m}_{\sigma_{f}\cap |CC({\cal F})|}(T^{*}M,or_{T^{*}M}) \to 
H^{2m}_{I}(T^{*}M,or_{T^{*}M})\to H^{2m}_{c}(T^{*}M,or_{T^{*}M}) \:.\]  

This intersection formula is an important and useful generalization of
\cite[thm.9.1, p.207]{Ka} and \cite[thm.9.5.6, p.386]{KS}. Note that we only assume,
that the intersection of $\sigma_{f}$ and the {\em support} $|CC({\cal F})|$
of the characteristic cycle of ${\cal F}$ is contained in a compact 
subanalytic set, and this is the reason why our formula is a result about
{\em Lagrangian cycles}. The results of loc.sit. are formulated under suitable
assumptions on the {\em micro-support} of ${\cal F}$, which is in general much
bigger than $|CC({\cal F})|$.

\begin{rem} The characteristic cycle $CC({\cal F})$ of ${\cal F}$
depends only on the {\em constructible function} $\alpha: x\mapsto \chi({\cal F}_{x})$ 
(compare with \cite[thm.9.7.11]{KS} and \cite[subsec.5.0.3]{Sch}).
Therefore our results can also be stated in terms of the associated
constructible function $\alpha$. In this context, one can only use the support
of the characteristic cycle $CC(\alpha)$ as an invariant of the constructible function
$\alpha$
(the micro-support $\mu supp({\cal F})$ depends on the sheaf complex
${\cal F}$, and {\em not} (!) only on the corresponding constructible function).
\end{rem}
  
Our proof of theorem \ref{thm:int1} is a modification of the proof of
\cite[thm.9.5.6]{KS}. We use the result of \cite[1.20, D.19]{vDrMi}
(due in the subanalytic context to Bierstone, Milman and Pawlucki),
that $K=\{g=0\}$ is the zero-set of a non-negative subanalytic $C^{2}$-function $g$.
Then there exists a relatively compact open neighborhood $U$ of $K$ in $M$, and a $\delta_{0}>0$ such that
\[g: \{g<\delta_{0}\}\cap U
\to [0,\delta_{0}[ \]
is proper (since $\{g=0\}=K$ is compact). After restriction to $U$ we can assume $M=U$.
Then
\begin{equation} \label{eq:Milnor1}
\chi\bigl(R\Gamma(K,R\Gamma_{\{f\geq 0\}}{\cal F})\bigr) = 
\chi\bigl(H^{*}(\{g\leq \delta\},M^{-}_{f},{\cal F})\bigr)
= \chi\bigl(H^{*}(\{g\leq \delta\},{\cal F})\bigr) -\chi\bigl(H^{*}(M^{-}_{f},{\cal F})\bigr)
\end{equation}
calculates for $0<\epsilon << \delta << \delta_{0}$ (i.e. for $\delta$ and $\epsilon$ sufficiently small,
with $\epsilon$ small compared to $\delta$) the Euler characteristic of the relative cohomology of a 
``tube'' $\{g\leq \delta\}$ modulo the {\em left Milnor fiber} 
\[M^{-}_{f}:=\{g\leq \delta, f=-\epsilon\} \quad (0<\epsilon << \delta << \delta_{0})\]
(compare \cite[sec.1.1]{Sch}).
In terms of the constructible function $\alpha$, this Euler characteristic can also be rewritten as
(compare \cite[sec.9.7]{KS} and \cite[sec.2.3]{Sch}):
\begin{equation} \label{eq:Milnor2}
\int_{\{g\leq \delta\}}\, \alpha \,d\chi - \int_{M^{-}_{f}}\, \alpha \,d\chi =
\int_{K}\, \alpha \,d\chi - \int_{M^{-}_{f}}\, \alpha \,d\chi \:.
\end{equation}

Here we use for a {\em compact} subanalytic subset $A\subset M$ the notation
\[\int_{A}\, \alpha \,d\chi := \int_{M}\, 1_{A}\cdot \alpha \,d\chi\:,
\quad \text{with} \quad 
\int_{M}\, (\cdot)  \,d\chi : \: CF_{c}(M,\bb{Z}) \to \bb{Z}\]
the {\em unique} $\bb{Z}$-linear map on the group $CF_{c}(M,\bb{Z})$ of $\bb{Z}$-valued
subanalytically constructible functions with {\em compact} support such that
\[\int_{M}\, 1_{B} \,d\chi = \chi\bigl(H^{*}(B,\bb{Q})\bigr)\]
for $B\subset M$ a compact subanalytic subset.\\

For the special case $f\equiv 0$ we get $\emptyset=M^{-}_{f}$ so that theorem~\ref{thm:int1}
implies the following generalization of the classical {\em Poincar\'{e}-Hopf index formula},
which corresponds to the case $\alpha=1_{M}$ for a $M$ a compact manifold
(with $\sigma_{0}=T^{*}_{M}M$ the zero-section of $T^{*}M$):

\begin{cor}[(global index formula)] \label{cor:globalindex}
Consider a constructible complex ${\cal F}$ as in theorem~\ref{thm:int1}
and denote by $\alpha$ the corresponding constructible function 
$x\mapsto \chi({\cal F}_{x})$. Assume ${\cal F}$ or $\alpha$ has a compact support.
Then 
\begin{equation} \label{eq:globalindex}
\chi\bigl(R\Gamma(M,{\cal F})\bigr) = \sharp\bigl(\; [T^{*}_{M}M]\cap [CC({\cal F})]\; \bigr) 
\quad \text{or} \quad 
\int_{M}\, \alpha \,d\chi = \sharp\bigl(\; [T^{*}_{M}M]\cap [CC(\alpha)]\; \bigr) \:.
\end{equation}
\end{cor}

This global index formula goes back to Dubson \cite[thm.2, p.115]{Du2} in the complex, and to
Kashiwara \cite[thm.8.1, p.205]{Ka} in the real context. Other presentations can be found in
\cite[thm.1.5, p.832]{Fu}, \cite[thm.9.1, p.382]{Gi}, \cite[cor.9.5.2, p.384]{KS}
and \cite[eq.(5.30), p.298]{Sch}.\\ 

Another important special case of our intersection formula (\ref{eq:int1})
is the case
$I=\{\omega\}$ given by a point $\omega\in T^{*}M$, with $x:=\pi(\omega)$ and $K:=\{x\}$. 
In this case we can use for $g$ the usual distance function (in suitable coordinates
$(M,x) \simeq (\bb{R}^{m},0)$). Then we get (cf. \cite[thm.9.5.6, p.386]{KS}): 
\begin{equation} \label{eq:int2}
\chi\bigl((R\Gamma_{\{f\geq 0\}}{\cal F})_{x}\bigr) = 
\sharp_{df_{x}}\bigl(\; [df(M)]\cap [CC({\cal F})]\; \bigr)
\end{equation}
and
\begin{equation} \label{eq:int2CF}
\alpha(x) - \int_{M^{-}_{f}}\, \alpha \,d\chi = 
\sharp_{df_{x}}\bigl(\; [df(M)]\cap [CC(\alpha)]\; \bigr)
\:,
\end{equation}
with $\sharp_{\omega}$ the corresponding {\em local intersection number}
and $M^{-}_{f}$ the {\em local left Milnor fiber} of $f$ in $x$.
We get in particular:
\begin{equation} \label{eq:int3}
\sigma_{f}\cap |CC({\cal F})| = \emptyset 
\:\:\Rightarrow \:\:
\chi\bigl((R\Gamma_{\{f\geq 0\}}{\cal F})_{x}\bigr) = 0\:.
\end{equation}

\begin{rem} In another paper, we will use the formula (\ref{eq:int3})
(together with a specialization result for Lagrangian cycles) for a
proof of the non-characteristic pullback formula \cite[prop.9.4.3, p.378]{KS}
under the {\em weaker} assumption, that the map is only non-characteristic with 
respect to the {\em support} 
$|CC({\cal F})|$ of the characteristic cycle of ${\cal F}$.
\end{rem}

As a further example, let us consider (locally) the case $f=r$ a (real analytic)
distance function to the point $x$, i.e. $r\geq 0$ with $\{x\}=\{r=0\}$.
By the {\em curve selection lemma} there exists an open neighborhood $U$ of $x$ in $M$,
with
\[dr(U)\cap |CC({\cal F})| \subset \{dr_{x}\} \quad \text{and} \quad
dr(U)\cap |CC(\alpha)| \subset \{dr_{x}\} \:.\]
This will be explained later on in terms of stratification theory
(and compare with \cite[prop.8.3.12, p.332]{KS}). So if we work on the manifold $U$,
then we can apply (\ref{eq:int2}) and (\ref{eq:int2CF}) with $\emptyset = M^{-}_{r}$:

\begin{cor}[(local index formula)] \label{cor:localindex}
\begin{equation} \label{eq:localindex}
\chi\bigl({\cal F}_{x}\bigr) = \sharp_{dr_{x}}\bigl(\; [dr(U)]\cap [CC({\cal F})]\; \bigr)
\quad \text{and} \quad 
\alpha(x)  = \sharp_{dr_{x}}\bigl(\; [dr(U)]\cap [CC(\alpha)]\; \bigr)
 \:.
\end{equation}
\end{cor}

So (\ref{eq:localindex}) describes an {inversion formula} for reconstructing the constructible
function $\alpha$ out of the Lagrangian cycles $[CC({\cal F})]$ and $[CC(\alpha)]$
(compare \cite[thm.2, p.115]{Du2}, \cite[thm.8.3, p.205]{Ka}, \cite[eq.(9.5.8), p.386]{KS}
and \cite[eq.(5.30), p.298]{Sch}).\\ 

Theorem \ref{thm:int1} also gives a purely {\em real} proof of the 
corresponding intersection formula in complex geometry:

\begin{cor} \label{cor:int2}
Let $M$ be an m-dimensional complex analytic manifold and 
$h: M\to \bb{C}$ a holomorphic function. 
Consider on $M$ a bounded complex analytically constructible 
complex
${\cal F}$ of sheaves of vector spaces, with finite 
dimensional stalks ${\cal F}_{x}$ ($x\in M$). Suppose that the intersection
of $\sigma_{h}:=\{(x,dh_{x})\in T^{*}M | x\in M\}$ and the support 
$|CC({\cal F})|$ of the characteristic cycle of ${\cal F}$ is 
contained in a
compact analytic subset $I\subset T^{*}M$, with $K:=\pi(I)\subset \{h=0\}$.
Then one has
\begin{equation} \label{eq:int4}
\chi\bigl(R\Gamma(K,\phi_{h}[-1]\;{\cal F})\bigr) = 
\sharp\bigl(\; [dh(M)]\cap [CC({\cal F})]\; \bigr) \:.
\end{equation}
\end{cor}

Here we use the notation $T^{*}M$ for the {\em holomorphic} cotangent bundle,
and $\phi_{h}$ is the {\em vanishing cycle functor} of Deligne.
The holomorphic section $dh$ of $T^{*}M$ corresponds under the natural
isomorphism $T^{*}M\simeq T^{*}(M^{R})$  to the section
$d(re(h))$ of $T^{*}(M^{R})$, with $T^{*}(M^{R})$ the {\em real}
cotangent bundle of the underlying real manifold $M^{R}$.
If we use the induced {\em complex orientation}
of $T^{*}(M^{R})$, then the class 
\[ [df(M^{R})]\in H^{2m}_{\sigma_{f}}(T^{*}M,or_{T^{*}(M^{R})/M^{R}})\simeq
H^{2m}_{\sigma_{h}}(T^{*}M,\bb{Z}) \]
for $f:=re(h)$ corresponds 
under Poincar\'{e} duality to the {\em fundamental class} in Borel-Moore homology
of the complex manifold $\sigma_{h}$.\\

Then the corollary follows from 
theorem \ref{thm:int1} and the isomorphism (cf.
\cite[lem.1.3.2, p.70]{Sch} and \cite[cor.1.1.1, p.32]{Sch}):
\begin{equation} \label{eq:vcyc}
R\Gamma(K,\phi_{h}[-1]\;{\cal F}) \simeq
R\Gamma(K,R\Gamma_{\{re(h)\geq 0\}}{\cal F})\:.
\end{equation}
Or in terms of the constructible function $\alpha$ :
\begin{equation} \label{eq:vcyc2}
\int_{M^{-}_{f}}\, \alpha \,d\chi = \int_{M_{h}}\, \alpha \,d\chi 
\end{equation}
and
\begin{equation} \label{eq:vcyc3}
\chi\bigl(R\Gamma(K,\phi_{h}[-1]\;{\cal F})\bigr) =
\int_{K}\, \alpha \,d\chi - \int_{M_{h}}\, \alpha \,d\chi \:,
\end{equation}
with 
\[M_{h}:=\{g\leq \delta, h=w\} \quad (0<|w|<< \delta << \delta_{0})\]
the {\em Milnor fiber} of the holomorphic function $h$.\\

This holomorphic intersection formula for the vanishing cycle functor is due 
to Dubson \cite[thm.1, p.183]{Du}, Ginsburg \cite[prop.7.7.1, p.378]{Gi}, 
L\^{e} \cite[thm.4.1.2, p.242]{Le2} and Sabbah \cite[thm.4.5, p.174]{Sa}. 
For a discussion of the history of this holomorphic intersection formula
we recommend the paper \cite{Le2}.

But most of these
references are in the language of (regular) {\em holonomic D-modules}, or {\em perverse
sheaves} (with respect to middle perversity). So the assumption on the 
intersection for a holomorphically constructible complex of sheaves 
corresponds to an assumption on the micro-support (\cite[thm.11.3.3, p.455]{KS}).
Only the result of Sabbah \cite[thm.4.5, rem.4.6, p.174]{Sa} is in terms of the 
underlying (complex analytic) Lagrangian cycle, which corresponds therefore
to an assumption on the intersection of $\sigma_{h}$ and the support of the
corresponding characteristic cycle! Similarly, only this reference \cite{Sa}, and in the
real subanalytic context also \cite[thm.9.1, p.207]{Ka} consider the case of an
intersection in a {\em compact} subset $I$. \\

The other references deal only with the case $I=\{\omega\}$ given by a point 
$\omega\in T^{*}M$. In this special case we 
get back a formula conjectured by Deligne (with $x:=\pi(\omega)$ and $K:=\{x\}$):
\begin{equation} \label{eq:int5}
\chi\bigl((\phi_{h}[-1]\;{\cal F})_{x}\bigr) =
\sharp_{dh_{x}}\bigl(\; [dh(M)]\cap [CC({\cal F})]\; \bigr)
\end{equation}
and
\begin{equation} \label{eq:int5CF}
\alpha(x) - \int_{M_{h}}\, \alpha \,d\chi = 
\sharp_{dh_{x}}\bigl(\; [dh(M)]\cap [CC(\alpha)]\; \bigr)
\:,
\end{equation}
with $M_{h}$ the {\em local Milnor fiber} of $h$ in $x$.
We get in particular:
\begin{equation} \label{eq:int6}
\sigma_{h}\cap |CC({\cal F})| = \emptyset 
\:\:\Rightarrow\:\:
\chi\bigl((\phi_{h}[-1]\;{\cal F})_{x}\bigr) = 0\:.
\end{equation}

We used this formula in our paper \cite{Sch2} for a short proof
of a formula of Brasselet, L\^{e} and Seade for the {\em Euler obstruction}
(\cite[thm.3.1]{BLS}).\\

As another application of the formula (\ref{eq:vcyc}) let us consider
the following example: ${\cal F}=k_{N}$, with $N\subset M$ a {\em closed} complex
analytic submanifold. Consider a holomorphic function germ $h: (M,x)\to
(\bb{C},0)$ such that $h|N$ has in $x\in N$ a {\em complex Morse critical  
point}. Then we get $CC(k_{N})=[T^{*}_{N}M]$, with the orientation convention
of \cite[chap. IX]{KS}. If we denote by 
\[[T^{*}_{N}M]_{c}\in H^{2m}_{T^{*}_{N}M}(T^{*}M,\bb{Z})\]
the class, which corrresponds under Poincar\'{e} duality to the fundamental
class of the complex manifold $T^{*}_{N}M$, then we get of course
\[\sharp_{dh_{x}}\bigl(\; [dh(M)]\cap [T^{*}_{N}M]_{c}\; \bigr)\;= 1\:.\]
But $re(h)|N$ has in $x$ a {\em real Morse critical point} of index 
$n:=dim_{C}(N)$. We therefore get by \cite[eq.(9.5.18), p.388]{KS} for $f:=re(h)$:
\[\sharp_{df_{x}}\bigl(\; [df(M)]\cap [CC(k_{N})]\; \bigr)\;= (-1)^{n}\:.\]
And this implies 
\[[T^{*}_{N}M]_{c} = (-1)^{n}\cdot [T^{*}_{N}M]  
\in H^{2m}_{T^{*}_{N}M}(T^{*}M,\bb{Z})\:.\]

So one has to be very careful about {\em orientation conventions} and the definition
of {\em characteristic cycles} that one uses. Here we follow the notations and conventions
of \cite{KS}. Note that there are many approaches to this subject,
often using quite different techniques and conventions. 
We recall in the next section a detailed comparison, which is worked out in \cite[subsec.5.0.3]{Sch}.\\

Let us only recall, that the theory of characteristic cycles has its origin in the theory
of {\em holonomic D-modules} and the {\em local index formula} of Kashiwara \cite[thm. on p.804]{Ka2}
(compare \cite[thm.6.3.1, p.127]{Ka2}, \cite[thm.2, p.574]{BDK} and \cite[thm.11.7, p.393]{Gi}).
This corresponds to corollary~\ref{cor:localindex}
for
\[{\cal F} = Rhom_{D_{M}}({\cal M},{\cal O}_{X}) =:sol({\cal M})\]
the {\em solution complex} of a holonomic $D_{M}$-module ${\cal M}$.

Later on, the theory of characteristic cycles was extended to constructible functions and sheaves
in the context of real geometry. First in the {\em subanalytic} context by Kashiwara \cite{Ka2,KS}, 
and independently also by Fu \cite{Fu}. A simple approach in the {\em semialgebraic} context is
sketched in \cite{GrM}, and the extension to {\em o-minimal structures} and 
{\em analytic geometric categories}
has been worked out in \cite{SchVi}. All these approaches in real geometry are based on a
suitable {\em Morse theory}, e.g. the {\em micro-local sheaf theory} of Kashiwara-Schapira \cite{KS},
or the {\em stratified Morse theory} of Goresky-MacPherson \cite{GM}.\\

For a detailed comparison and translation of these two different theories
in the framework of {\em Morse theory for constructible sheaves},
including a {\em geometric} introduction to characteristic cycles of constructible functions and sheaves,
we refer to our book \cite[chapter 5]{Sch}. This language will be used in the 
proof of theorem \ref{thm:int1}.
Moreover, we use the following important result about the behaviour
of the {\em support}  
$|CC({\cal F})|$ of the characteristic cycle of ${\cal F}$ under a
suitable intersection:

\begin{thm} \label{thm:int2}
Let $M$ be a real analytic manifold and $g: M\to \bb{R}$ a subanalytic
$C^{2}$-function. Consider on $M$ a bounded subanalytically constructible 
complex
${\cal F}$ of sheaves of vector spaces, with finite 
dimensional stalks ${\cal F}_{x}$ ($x\in M$). 
Assume ${\cal F}$ is constructible with respect to a subanalytic
Whitney b-regular stratification (i.e. the cohomology sheaves of
${\cal F}$ are locally constant on the subanalytic strata $S$) such that the 
set of nondegenerate covectors (with respect to this stratification) is 
dense in all fibers of the projection $T^{*}_{S}M \to S$ (for all strata 
$S$). Let $\delta$ be a regular value of $g$ such that $\{g=\delta\}$ is 
{\em transversal} to all strata $S$. Then one has for the open inclusion
$j: \{g<\delta\}\to M$ the following estimate for the support 
$|CC(Rj_{*}j^{*}{\cal F})|$ of the characteristic cycle of 
$Rj_{*}j^{*}{\cal F}$:
\begin{equation} \label{eq:est*}
\begin{split}
|CC(&Rj_{*}j^{*}{\cal F})| \cap \pi^{-1}(\{g=\delta\}) \subset
\\
\{\;\omega+\lambda\cdot dg_{x}\;\in T^{*}M \;|&\; \lambda \leq 0,\;
\pi(\omega)=x,\; g(x)=\delta,\; \omega\; \in |CC({\cal F})|\;\}\:.
\end{split} 
\end{equation}
\end{thm}

\begin{rem} The assumption on our stratification is for example satisfied 
for a subanalytic {\em $\mu$-stratification} in the sense
of \cite[def.8.3.19, p.334]{KS} (see \cite[cor.8.3.24, p.336]{KS}). 
This {\em $\mu$-condition} is by \cite{Tr} equivalent to the {\em w-condition} of Verdier,
and such a {\em w-regular} stratification can always be used in the context of
``geometric categories'' as in the last section of this paper (cf. \cite{Loi}).
But it is not known (at least to the author), if the {\em w-condition} implies
the assumption on our stratification (used in the theorem) in this more general context.
\end{rem}

Kashiwara-Schapira used in their proof of
\cite[thm.9.5.6, p.386]{KS} a similar result in terms of the {\em micro-support} of the
corresponding constructible complex of sheaves 
(cf. \cite[prop.8.4.1, p.338]{KS} and \cite[prop.5.4.8.(a), p.233]{KS}).

We will give a proof based on our results on {\em Morse theory
for constructible sheaves} \cite[chapter~5]{Sch}, 
especially on an explicit
result about suitable ``stratified spaces with boundary'' \cite[thm.5.0.2, p.286]{Sch}. 
This approach works also, if one uses the supports of
the corresponding characteristic cycles. It implies a similar
estimate for 
\[|CC(Rj_{!}j^{*}{\cal F})| \cap \pi^{-1}(\{g=\delta\}) \:,\]
(with $\lambda\geq 0$), from which we deduce (\ref{eq:est*}) by {\em duality}.


\section{Characteristic cycles}

In this section we explain the construction of the characteristic cycles.
As an application we give a proof of theorem \ref{thm:int2} in a much
more general context, which applies especially to complexes of sheaves
constructible with respect to suitable ``geometric categories''.\\

Let $X$ be a closed subset of a smooth manifold $M$.
In this paper, a smooth manifold $M$ has by definition a countable topology.
We also assume (for simplicity), that the dimension of the connected
components of $M$ is bounded.
A {\em stratification} of $X$ is a filtration $X_{\bullet}$ of $X$:
\[\emptyset=X_{-1}
\subset X_{0}\subset \dots \subset X_{n}=X\] 
by closed subsets such that $X^{i}:=X_{i}\backslash X_{i-1}$ ($i=0,\dots,n$) is 
a smooth submanifold of $M$. The connected components of the $X^{i}$ are by 
definition the {\em strata} of this stratification.

Here we fix a degree $k = 1,\dots,\infty,\omega$ of smoothness
(with $\omega =$ real analytic).
(Sub)manifold and smooth map or function always means a 
$C^{k}$-(sub)manifold and $C^{k}$-smooth map or function.
We assume that the stratification $X_{\bullet}$ is {\em Whitney b-regular}:
\begin{enumerate}
\item[]
If $x_{n}\in X^{j}$ (for $i<j$) and $y_{n}\in X^{i}$ are sequences converging 
to $x\in X^{i}$ such that the tangent planes $T_{x_{n}}X^{j}$ converge to some 
limiting plane $\tau$, and the secant lines $l_{n}=\overline{x_{i},y_{i}}$ 
(with respect to some local coordinates) converge to some limiting line $l$, 
then $l \subset\tau$.
\end{enumerate}

Consider on $X$ a bounded complex
${\cal F}$ of sheaves of vector spaces (over a base-field $k$), with finite 
dimensional stalks ${\cal F}_{x}$ ($x\in X$). 
We assume that ${\cal F}$ is {\em constructible} with respect to the stratification 
$X_{\bullet}$, i.e. the cohomology sheaves of
${\cal F}$ are {\em locally constant} on all $X^{i}$.\\

From now on we assume $X=M$ and $k\geq 2$. 
Then we have by \cite[cor.4.0.3, p.216]{Sch} and \cite[prop.4.1.2, p.227]{Sch} the following 
estimate for
the {\em micro-support} $\mu supp({\cal F})$ of ${\cal F}$ in the sense of
\cite[def.5.1.2, p.221]{KS}:
\begin{equation} \label{eq:est.musupp}
\mu supp({\cal F}) \subset \Lambda:=\bigcup_{i} \: T_{X^{i}}^{*}X 
\subset T^{*}X \:.
\end{equation}
Here $T_{S}^{*}X$ denotes the conormal bundle of a locally closed submanifold 
$S$ of $X$. By Whitney regularity, $\Lambda$ is a closed
subset of $T^{*}X$ (this is in fact equivalent to {\em a-regularity}).\\

Moreover, we have by \cite[prop.4.0.2, p.219]{Sch} that 
${\cal F}$ is {\em cohomologically constructible} in the sense of
\cite[def.3.4.1, p.158]{KS}. Therefore we can define as in \cite[p.377]{KS}
the following chain of morphisms:
\[RHom({\cal F},{\cal F})\simeq R\pi_{*}\mu hom({\cal F},{\cal F})
\simeq R\pi_{*}R\Gamma_{\Lambda}\mu hom({\cal F},{\cal F})\simeq\]
\[R\pi_{*}R\Gamma_{\Lambda}\mu_{\triangle}({\cal F}\boxtimes D({\cal F}))
\to R\pi_{*}R\Gamma_{\Lambda}\mu_{\triangle}\delta_{*}
({\cal F}\otimes D({\cal F}))\to\]
\[R\pi_{*}R\Gamma_{\Lambda}\mu_{\triangle}\delta_{*}\omega_{X} \simeq
R\pi_{*}R\Gamma_{\Lambda}(\pi^{-1}\omega_{X}) \;,\]
with $\pi: T^{*}X\to X$ the projection, $\delta: X\to X\times X$ the
diagonal embedding and $\omega_{X}$ the dualizing
complex on $X$. Since $X$ is a manifold, one has by \cite[prop.3.3.2, p.152]{KS}:
\[\omega_{X} \simeq or_{X}[d]\;, \quad \text{with $d:=dim(X)$.}\]
Here we assume that $X$ is {\em pure-dimensional} (or one should interpret
$d:=dim(X)$ as a locally constant function).\\ 

For the definition of the above chain of morphisms we use the following
properties:
\begin{enumerate}
\item $RHom({\cal F},{\cal F})\simeq R\pi_{*}\mu hom({\cal F},{\cal F})$
(\cite[prop.4.4.2(i), p.202]{KS}).
\item $\mu hom({\cal F},{\cal F})$ has its support in $\mu supp({\cal F})$
(\cite[cor.5.4.10.(ii), p.234]{KS}), and therefore by (\ref{eq:est.musupp}) also in
$\Lambda$.
\item $\mu hom({\cal F},{\cal F})\simeq
\mu_{\triangle}({\cal F}\boxtimes D({\cal F}))$ (\cite[def.4.4.1.(iii),
p.202]{KS}
and \cite[prop.3.4.4, p.159]{KS}). Here we can apply \cite[prop.3.4.4]{KS}, since
${\cal F}$ is cohomologically constructible.
\item We finally use the natural morphisms
${\cal F}\boxtimes D({\cal F})\to 
\delta_{*}\delta^{*}({\cal F}\boxtimes D({\cal F}))\simeq
\delta_{*}({\cal F}\otimes D({\cal F}))$,
and ${\cal F}\otimes D({\cal F})\to
\omega_{X}$, together with $\mu_{\triangle}\delta_{*}\omega_{X} \simeq
\pi^{-1}\omega_{X}$.
\end{enumerate}

\begin{defi} \label{def:CC}
The image of $id\in$ Hom$({\cal F},{\cal F})\simeq H^{0}(RHom({\cal F},{\cal F}))$ in 
$H^{0}_{\Lambda}(T^{*}X,\pi^{-1}\omega_{X})
\simeq H^{d}_{\Lambda}(T^{*}X,\pi^{-1}or_{X})$ is called the {\em characteristic
cycle} $CC({\cal F})$ of ${\cal F}$.
\end{defi}

\begin{rem} This definition extends \cite[def.9.4.1, p.377]{KS} to our context.
The same construction works for a Whitney {\em a-regular stratification}
$X_{\bullet}$, if we have the estimate (\ref{eq:est.musupp}), and the 
property that ${\cal F}$ is cohomologically constructible, e.g. $X_{\bullet}$ is 
(locally) {\em C-regular} in the sense of Bekka \cite{Be} (compare with
\cite[chapter~4]{Sch}).
\end{rem}

We now explain the calculation of $CC({\cal F})$ in terms of our (stratified)
{\em Morse theory for constructible sheaves} \cite[chapter~5]{Sch}.\\

Consider the set of {\em nondegenerate} covectors with respect to our
Whitney b-regular stratification $X_{\bullet}$ 
(see \cite[def.5.1.2, p.308]{Sch} and \cite[def.1.8, p.44]{GM}):
\begin{equation} \label{eq:regcov}
\Lambda':= \bigcup_{i} \: \Bigl(T_{X^{i}}^{*}X \backslash \bigcup_{i\neq j} \:
cl(T_{X^{j}}^{*}X) \Bigr) \:.
\end{equation}

Choose a smooth function germ $f: (X,x)\to (\bb{R},0)$, with $x\in X^{s}$
and $df_{x}\in \Lambda'$. Take a {\em normal slice} $N$ at $x$,
i.e. a locally closed submanifold $N$ of $X$,
with $N\cap X^{s}=\{x\}$ such that $N$ intersects $X^{s}$ transversally in
$x$ (cf. \cite[def.5.0.2(2), p.277]{Sch} and \cite[def.1.4, p.41]{GM}). 
Then the isomorphism-class of
\begin{equation} \label{eq:NMD}
NMD({\cal F},f,x):=
\bigl(R\Gamma_{\{f|N\geq 0\}}({\cal F}|N)\bigr)_{x}
\end{equation} 
is the (sheaf theoretic) {\em ``normal Morse datum''} of
$f$ in $x$ with respect to ${\cal F}$, i.e. it is the cohomological counterpart
of the corresponding {\em normal Morse data} of Goresky-MacPherson \cite[def.3.6.1, p.65]{GM}.
By \cite[thm.5.0.1(2), p.278]{Sch},
this isomorphism-class depends {\em only} on $df_{x}$ and is {\em locally constant}
on $\Lambda'$.
By \cite[prop.5.0.2, p.279]{Sch} we have a distinguished triangle
\begin{equation} \begin{CD}
R\Gamma(l^{-}_{X},{\cal F})[-1] @>>> 
NMD({\cal F},f,x) @>>>
{\cal F}_{x}  @> [1] >> \: ,
\end{CD} \end{equation}
with $l^{-}_{X}$ the {\em ``lower halflink''} of $f$ in $x$ (\cite[def.3.9.1, p.66]{GM}).
In the notation of this paper, this $l^{-}_{X}$ is just the intersection of $X$ with a
{\em local left Milnor fiber} of $f|N$
in $x$. But ${\cal F}|l^{-}_{X}$ is constructible with respect to an induced
Whitney b-regular stratification of $l^{-}_{X}$. Then
\[R\Gamma(l^{-}_{X},{\cal F})\]
has {\em finite
dimensional} cohomology, since $l^{-}_{X}$ is compact (cf. 
\cite[rem.4.2.2, p.245]{Sch} and \cite[thm.3.5, p.70]{Intersection}).
By the above distinguished triangle, this is also true for $NMD({\cal F},f,x)$, with
\begin{equation} \label{eq:multi}
\chi\bigl(NMD({\cal F},f,x)\bigr) = \chi\bigl({\cal F}_{x}) - 
\chi\bigl(R\Gamma(l^{-}_{X},{\cal F})\bigr) \:.
\end{equation}

So we can associate to a connected component $\Lambda'_{j}$ of 
$\Lambda'$ the integer (for $df_{x}\in \Lambda'_{j}$):
\begin{equation} \label{eq:chiNMD}
m_{j}:=\chi(NMD({\cal F},f,x)):=
\chi\bigl(\bigl(R\Gamma_{\{f|N\geq 0\}}({\cal F}|N)\bigr)_{x}\bigr)\:.
\end{equation}

If we work in the context of ``geometric categories'' as in the next section, with
$f$ a {\em definable} function germ, then this multiplicity $m_{j}$ 
can also be expressed in terms of the {\em constructible function}
$\alpha: x\mapsto \chi({\cal F}_{x})$ (compare \cite[eq.(5.16), p.290]{Sch}):
\begin{equation} \label{eq:chiNMD2}
m_{j}:=\chi(NMD(\alpha,f,x)):=
\alpha(x) - \int_{l^{-}_{X}}\, \alpha \,d\chi \:.
\end{equation}

Take now an open subset $\Omega$ in $T^{*}X$, with $\Omega\cap \Lambda
= \Lambda'$ (e.g. the complement of the set 
$\Upsilon:=\Lambda\backslash \Lambda'$ of {\em degenerate} covectors).
Then one gets for the image of $CC({\cal F})$ under the natural map
$H^{d}_{\Lambda}(T^{*}X,\pi^{-1}or_{X})\to
H^{d}_{\Lambda'}(\Omega,\pi^{-1}or_{X})$ the formula
\begin{equation} \label{eq:CC}
im(\;CC({\cal F})\;) = \prod \: m_{j}\cdot [\Lambda'_{j}]\:.
\end{equation}
Here $[\Lambda'_{j}]$ is defined as in \cite[chap.9.4]{KS}, i.e.
it is the image of the class $[T_{X^{s}}^{*}X] \in 
H^{d}_{T_{X^{s}}^{*}X}(\Omega,\pi^{-1}or_{X})$ of loc.sit.
(with $\pi(\Lambda'_{j})\subset
X^{s}$) under the projection
\[H^{d}_{T_{X^{s}}^{*}X}(\Omega,\pi^{-1}or_{X})\simeq
\prod \: H^{d}_{\Lambda'_{k}}(\Omega,\pi^{-1}or_{X}) 
\to H^{d}_{\Lambda'_{j}}(\Omega,\pi^{-1}or_{X})\:,\]
where the product is over all $k$ with $\pi(\Lambda'_{k})\subset
X^{s}$.\\

This follows as in \cite[p.382]{KS} from our identification \cite[eq.(5.52), p.318]{Sch} of the 
{\em ``normal Morse datum''} $NMD({\cal F},f,x)$ of
$f$ in $x$ with the {\em ``local type of ${\cal F}$ in $df_{x}$''} in the sense of
Kashiwara-Schapira (compare \cite[6.6.1(ii), p.274]{KS} and \cite[def.7.5.4, p.311]{KS}).\\

Assume $H^{d}_{\Upsilon}(T^{*}X,\pi^{-1}or_{X}) = 0$ so that the above map of
local cohomology groups is injective. Then the characteristic cycle
$CC({\cal F})$ is uniquely determined by the Euler characteristics 
$m_{j}$ of the ``normal Morse data'' of ${\cal F}$.
This applies for example, if $\Upsilon$ has a stratification with all strata
of dimension $< d=dim(X)$ (e.g. in the context of ``geometric categories''
as in the next section). It would also be enough
(compare with the proof of \cite[prop.9.2.2.(i), p.367]{KS}), that $\Upsilon$ has 
locally this property (i.e. each point in $\Upsilon$ has an open
neighborhood $U$ in $T^{*}X$ such that $U\cap \Upsilon$ has such a 
stratification).

\begin{rem} By (\ref{eq:multi}) and (\ref{eq:chiNMD2}) we get in this case
an easy {\em geometric} description of the characteristic cycle $CC(\cdot)$,
or more precisely, of $im\bigl(\;CC(\cdot)\;)$. And one should ask if one 
can use this description as a {\em definition} of the characteristic cycle
$CC(\cdot)$. The main problem is then to show that this is a {\em cycle}
coming from 
\[H^{d}_{\Lambda}(T^{*}X,\pi^{-1}or_{X})\hookrightarrow
H^{d}_{\Lambda'}(\Omega,\pi^{-1}or_{X})  \:.\]
The sophisticated definition \ref{def:CC} seems to be the only one working
for a general {\em Whitney b-regular} stratification. But in the {\em subanalytic} context,
other more {\em geometric} approaches to this question are due to
Kashiwara \cite[thm.4.1, p.199]{Ka} and Fu (compare with \cite[def.4.1, p.856]{Fu} 
and \cite[thm.4.7, p.865]{Fu}), where \cite{Fu} uses the language of ``geometric measure theory''.
A ``translation'' of the last approach into ``geometric categories''
and a {\em specialization result in homology} is worked out in \cite[thm.5.0.3, p.296]{Sch}
and \cite[subsec.5.2.2]{Sch}.
\end{rem}

Suppose that $\Upsilon$ has (locally) such a stratification,
and that the set $\Lambda'$ of nondegenerate covectors is {\em dense} in $\Lambda$.
Then one gets the following ``explicit'' description of the {\em support}
$|CC({\cal F})|$ of the characteristic cycle of ${\cal F}$:
\begin{equation} \label{eq:suppCC}
|CC({\cal F})| = \bigcup_{m_{j}\neq 0} \:cl(\Lambda'_{j})\:.
\end{equation}

\begin{rem} \label{rem:musupp}
The micro-support $\mu supp({\cal F})$ contains by definition
the closure $cl(\Lambda'_{j})$ of all $\Lambda'_{j}$ such that the 
corresponding ``normal Morse datum'' 
\[NMD({\cal F},f,x) \simeq \bigl(R\Gamma_{\{f|N\geq 0\}}({\cal F}|N)\bigr)_{x}
\quad (df_{x}\in
\Lambda'_{j})\] 
is not isomorphic to $0$. Moreover, in the context
of ``geometric categories'', $\mu supp({\cal F})$ is {\em exactly} given as the
union of these $cl(\Lambda'_{j})$ (\cite[prop.5.0.1, p.279]{Sch}).
We see in particular, that in general 
$|CC({\cal F})|$ is much smaller than $\mu supp({\cal F})$.
\end{rem}

Suppose that the set $\Lambda'$ of nondegenerate covectors is dense in 
$\Lambda$, and that $\Upsilon$ has (locally) a stratification with all strata
of dimension $< d=dim(X)$ such that the frontier 
$\partial \Lambda'_{j}:=cl(\Lambda'_{j}) \backslash \Lambda'_{j}$
of $\Lambda'_{j}$ is a union of strata of $\Upsilon$ (for all $j$).
Then all morphisms in the following commutative diagram are injective:
\begin{displaymath} \begin{CD}
H^{d}_{|CC({\cal F})|}(T^{*}X,\pi^{-1}or_{X}) @>>>
H^{d}_{\Omega\cap |CC({\cal F})|}(\Omega,\pi^{-1}or_{X}) \\
@VVV  @VVV\\
H^{d}_{\Lambda}(T^{*}X,\pi^{-1}or_{X}) @>>>
H^{d}_{\Lambda'}(\Omega,\pi^{-1}or_{X})\:.
\end{CD} \end{displaymath}
 
Moreover, there exists a {\em unique} class $[CC({\cal F})]\in
H^{d}_{|CC({\cal F})|}(T^{*}X,\pi^{-1}or_{X})$, with image
$CC({\cal F}) \in  H^{d}_{\Lambda}(T^{*}X,\pi^{-1}or_{X})$, and image
\[\prod_{m_{j}\neq 0}\: m_{j}\cdot [\Lambda'_{j}]
\in H^{d}_{\Omega\cap |CC({\cal F})|}(\Omega,\pi^{-1}or_{X})\:.\]
This is the cohomology class that we used in the introduction (see the
next section, and compare with \cite{SchVi}). For the convenience of the reader,
we also compare our definition of the characteristic cycle $CC(\cdot)$  
with the other conventions used in the literature (compare with
\cite[subsec.5.0.3.]{Sch} for the details).\\

For simplicity, we assume that $M$ is {\em oriented}
(as in most of the following references), i.e. an isomorphism
$or_{M}\simeq \bb{Z}_{M}$ has been chosen.
Moreover, we {\em orient} $T^{*}M$ by using first
this orientation of the base $M$, and then the induced (real dual) orientation
of the fibers $T^{*}_{x}M=hom_{R}(T_{x}M,\bb{R})$ (as in \cite{Fu, GrM, SchVi, Sch}).
This differs by the factor $(-1)^{m(m+1)/2}$ from the {\em symplectic}
orientation of $T^{*}M$ (used in \cite{Ka}). In particular, for $M$ a complex analytic
manifold of complex dimension $m$, this agrees only up to the factor $(-1)^{m}$
with the {complex} orientation of $T^{*}M$! Let $[T^{*}M]$ be the {\em fundamental class}
of this oriented manifold in {\em Borel-Moore homology}, and denote by
$a: T^{*}M\to T^{*}M$ the {\em antipodal map} (i.e. multiplication by $-1$ on the fibers).

\begin{itemize}
\item[$\bullet$] By definition $CC({\cal F})$ corresponds to the characteristic cycle used in \cite{KS}.
\item[$\bullet$] $CC({\cal F})\cap (-1)^{m}\cdot [T^{*}M]$
corresponds to the characteristic cycle used 
in \cite{Ka, SchVi, Sch}.
\item[$\bullet$] $a_{*}\bigl(\; CC(\alpha)\cap (-1)^{m}\cdot [T^{*}M]\; \bigr)$
corresponds to the 
characteristic cycle used in \cite{Fu, GrM}.
\end{itemize}
Assume $M$ is a complex analytic manifold of complex dimension $m$,
with $[T^{*}M]_{c}$ the {\em fundamental class} of the complex manifold
$T^{*}M$. Let ${\cal F}$ or $\alpha$ be complex analytically constructible.
\begin{itemize}
\item[$\bullet$] $CC({\cal F})\cap (-1)^{m}\cdot [T^{*}M]_{c}$
corresponds to the characteristic cycle used 
in \cite{Du2, Du, Le2}.
\item[$\bullet$] $CC(\alpha)\cap [T^{*}M]_{c}$
corresponds to the characteristic cycle used 
in \cite{Sa}.
\item[$\bullet$] $CC({\cal F})\cap (-1)^{m}\cdot [T^{*}M]_{c}$, with 
${\cal F}:=sol({\cal M})$ or  
${\cal F}:=DR({\cal M}):=
Rhom_{D_{M}}({\cal O}_{X},{\cal M})$,
corresponds to the characteristic cycle of the {\em holonomic D-module}
${\cal M}$ used in \cite{BDK, Du2, Le2}.
\item[$\bullet$] $CC({\cal F})\cap [T^{*}M]_{c}$, with 
${\cal F}:=DR({\cal M}):=
Rhom_{D_{M}}({\cal O}_{X},{\cal M})[m]$,
corresponds to the characteristic cycle of the {\em holonomic D-module}
${\cal M}$ used in \cite{Gi}.
\end{itemize}

Now we prove the first main result of this section:

\begin{thm} \label{thm:int3}
Let $g: X\to \bb{R}$ be a $C^{2}$-function on the smooth manifold $X$. 
Consider a Whitney $b$-regular stratification $X_{\bullet}$ of $X$  
such that the 
set of nondegenerate covectors (with respect to this stratification) is
dense in all fibers of the projection $T^{*}_{S}X \to S$ (for all strata 
$S$). Let $\delta$ be a regular value of $g$ such that $\{g=\delta\}$ is 
transversal to all strata $S$. Assume that the set of degenerate covectors
of $X_{\bullet}$ and also of the induced Whitney $b$-regular stratification
of $\{g\leq \delta\}$ has (locally) a stratification with all strata
of dimension $< d=dim(X)$.

Consider on $X$ a bounded complex
${\cal F}$ of sheaves of vector spaces (over a base-field $k$), with finite 
dimensional stalks ${\cal F}_{x}$ ($x\in X$). 
Assume ${\cal F}$ is constructible with respect $X_{\bullet}$.
Then one has for the open inclusion
$j: \{g<\delta\}\to X$ the following estimate for the support of the 
characteristic cycles of $Rj_{!}j^{*}{\cal F}$ and 
$Rj_{*}j^{*}{\cal F}$:
\begin{equation} \label{eq:est!2}
\begin{split}
|CC(&Rj_{!}j^{*}{\cal F})| \cap \pi^{-1}(\{g=\delta\}) \subset \\
\{\;\omega+\lambda\cdot dg_{x}\;\in T^{*}M \;|&\;  \lambda \geq 0,\;
\pi(\omega)=x,\; g(x)=\delta,\;\omega\; \in |CC({\cal F})|\;\}\:, 
\end{split}
\end{equation}
\begin{equation} \label{eq:est*2}
\begin{split}
|CC(&Rj_{*}j^{*}{\cal F})| \cap \pi^{-1}(\{g=\delta\}) \subset \\
\{\;\omega+\lambda\cdot dg_{x}\;\in T^{*}M \;|&\;  \lambda \leq 0,\;
\pi(\omega)=x,\; g(x)=\delta,\;\omega\; \in |CC({\cal F})|\;\}\:.
\end{split}
\end{equation}
\end{thm}

\begin{proof} 1.) By transversality, $\{g\leq \delta\}$ gets an induced
Whitney $b$-regular stratification with strata $S'$ of the form 
$\{g<\delta\}\cap S,\;\{g=\delta\}\cap S$ (for $S$ a stratum of 
$X_{\bullet}$). 
The assumptions imply, that also the set of nondegenerate 
covectors with respect to this induced stratification of
$\{g\leq \delta\}$ is 
{\em dense} in all fibers of the projection $T^{*}_{S'}X \to S'$ (for all strata 
$S'$, compare with the proof of \cite[thm.5.0.2, p.286]{Sch}).\\

2.) $Rj_{!}j^{*}{\cal F}$ and $Rj_{*}j^{*}{\cal F}$ have also finite
dimensional stalks, and are constructible with respect to the induced 
stratification of $\{g\leq \delta\}$ (\cite[prop.4.0.2, p.219]{Sch}).
More precisely, they have their support in $\{g\leq \delta\}$,
and their restrictions to $\{g\leq \delta\}$ are constructible with respect 
to the induced stratification. 
Therefore the characteristic cycles of these complexes of sheaves
are defined, and we can apply the description (\ref{eq:suppCC}) for their
support.\\

3.) First we consider $Rj_{!}j^{*}{\cal F}$. Let 
\[\omega\;\in T^{*}_{S'}X \subset 
|CC(Rj_{!}j^{*}{\cal F})| \cap \pi^{-1}(\{g=\delta\})\]
be given, with $\omega$ {\em nondegenerate} with
respect to the induced stratification of $\{g\leq \delta\}$, and 
$x:=\pi(\omega)$ a point in a ``boundary'' stratum $S'= S\cap\{g= \delta\}$.
We can approximate $\omega$
by a covector $\lambda\cdot dg_{x} + \omega'\;\in T^{*}_{S'}X$ 
(with $\lambda \neq 0$ and $x=\pi(\omega')$)
such that $\omega'\;\in T^{*}_{S}X$ is {\em nondegenerate} 
with respect to $X_{\bullet}$, and the following
properties hold (see for instance \cite[thm.5.0.2, p.286]{Sch}):
\begin{enumerate}
\item[a.] $\lambda < 0$ implies that the ``normal Morse datum''
with respect to $Rj_{!}j^{*}{\cal F}$ at the covector $\omega$ is isomorphic
to $0$.
\item[b.] $\lambda > 0$ implies that the ``normal Morse datum''
with respect to $Rj_{!}j^{*}{\cal F}$ at the covector $\omega$ is isomorphic
to the ``normal Morse datum''
with respect to ${\cal F}$ at the covector $\omega'$, up to a shift by $[-1]$.
\end{enumerate}
Compare with \cite[p.286]{Sch} for the details (and note, that we used in
loc.cit. the function $\delta-g$ for the description of $\{g\leq \delta\}
=\{\delta-g\geq 0\}$). These results can also be deduced from 
\cite[prop.6.1.9(ii), p.256]{KS} and \cite[prop.7.5.10, p.314]{KS}.
So if the first ``normal Morse datum'' (at $\omega$) has a non-zero
Euler characteristic, then the same is true for the second 
``normal Morse datum'' (at $\omega'$).
By the description (\ref{eq:suppCC}), this implies our claim for 
$Rj_{!}j^{*}{\cal F}$ (since $Rj_{!}j^{*}{\cal F}$ and ${\cal F}$ have
the same ``normal Morse data'' for the strata $S'$ of the form
$\{g<\delta\}\cap S$).\\

4.) The claim for $Rj_{*}j^{*}{\cal F}$ follows from the case before
by the {\em duality} isomorphism
(with $D(\cdot)$ the duality functor as in \cite[def.3.1.16(ii), p.148]{KS}): 
\[D(Rj_{!}j^{*}{\cal F})\simeq Rj_{*}j^{*}(D({\cal F}))\:,\]
because $|CC(D({\cal G}))|$ ( for ${\cal G}={\cal F},\;Rj_{!}j^{*}{\cal F})$ ) 
is equal to the image of $|CC({\cal G})|$ under
the antipodal map $a: T^{*}X\to T^{*}X$ (i.e. multiplication by $-1$ on the 
fibers of $\pi$). This follows from the 
description (\ref{eq:suppCC}) together with \cite[eq.(5.25), p.296]{Sch}
(compare also with \cite[prop.9.4.4, p.380]{KS}).
\end{proof}

\begin{rem} Suppose that we are in the context of ``geometric categories''
(as in the next section). Then we can use \cite[prop.5.0.1, p.279]{Sch}
for the description of the {\em micro-support} (i.e. the description of remark
\ref{rem:musupp}), and the above proof gives also the corresponding
estimate in terms of the micro-support (instead of the support of the
characteristic cycles):
\[\mu supp(Rj_{!}j^{*}{\cal F}) \cap \pi^{-1}(\{g=\delta\}) \subset\]
\[\{\;\omega+\lambda\cdot dg_{x}\;\in T^{*}M \;|\;  \lambda \geq 0,\;
\pi(\omega)=x,\; g(x)=\delta,\;\omega\; \in \mu supp({\cal F})\;\}\;,\] 
\[\mu supp(Rj_{*}j^{*}{\cal F}) \cap \pi^{-1}(\{g=\delta\}) \subset\]
\[\{\;\omega+\lambda\cdot dg_{x}\;\in T^{*}M \;|\;  \lambda \leq 0,\;
\pi(\omega)=x,\; g(x)=\delta,\;\omega\; \in \mu supp({\cal F})\;\}\;.\]
In this way one can get in particular a proof of \cite[thm.5.3, p.293]{Ma} 
in terms of our (stratified) {\em Morse theory for constructible sheaves} \cite[chapter 5]{Sch},
without the use of the general {\em microlocal theory} of Kashiwara-Schapira. 
\end{rem}
 
For the proof of the definable counterpart of theorem \ref{thm:int1}, 
we need in the next section also (the second part of) the following result 
(compare with \cite[thm.4.2, thm.4.3, p.199]{Ka} and \cite[thm.9.5.3, cor.9.5.4, p.385]{KS}):

\begin{thm} \label{thm:int4}
Let $f: X\to [a,b[\;\subset \bb{R}$ be a $C^{2}$-function on the smooth 
manifold $X$ (with $b\in \bb{R}\cup \{\infty\}$). 
Consider a Whitney $b$-regular stratification $X_{\bullet}$ of $X$  
such that the set $\Lambda'$ of nondegenerate covectors (with respect to this 
stratification) is dense in $\Lambda$.
Assume that the set $\Upsilon:=\Lambda\backslash \Lambda'$ of degenerate 
covectors of $X_{\bullet}$ has (locally) a stratification with all strata
of dimension $< d=dim(X)$ such that the frontier 
$\partial \Lambda'_{j}:=cl(\Lambda'_{j}) \backslash \Lambda'_{j}$
is a union of strata of $\Upsilon$ (for all connected
components $\Lambda'_{j}$ of $\Lambda'$).

Let ${\cal F}$ be a bounded complex of sheaves
of vector spaces on $X$, with finite dimensional stalks ${\cal F}_{x}$ 
($x\in X$). Assume that ${\cal F}$ is constructible with respect to  
$X_{\bullet}$, and that $f|supp({\cal F})$ is proper.

\begin{enumerate}
\item Suppose that $df(X)\cap \mu supp({\cal F})$ is compact.
Then $R\Gamma(X,{\cal F})$ is a cohomologically bounded complex with finite 
dimensional cohomology, and
\[\chi(R\Gamma(X,{\cal F})) = \sharp(\:[df(X)]\cap [CC({\cal F})]\:)\:.\]
\item Suppose that $-df(X)\cap \mu supp({\cal F})$ is compact.
Then $R\Gamma_{c}(X,{\cal F})$ is a cohomologically bounded complex with 
finite dimensional cohomology, and
\[\chi(R\Gamma_{c}(X,{\cal F})) = \sharp(\:[-df(X)]\cap [CC({\cal F})]\:)\:.\]
\end{enumerate} \end{thm}

Here we use the definition of the {\em intersection product} of
the introduction, with $df(X):=\sigma_{f}:=\{(x,df_{x})\in T^{*}X | x\in X\}$.
Note, that we have the inclusion $|CC({\cal F})|\subset \mu supp({\cal F})$.

\begin{proof} The proof is by {\em Morse theory}, and is almost the same as the 
proof of \cite[thm.5.0.4, p.297]{Sch} given in \cite[p.327-328]{Sch} (compare also with \cite{Ka}
and the proof of \cite[thm.10.3.8, p.429]{KS}).

By approximation (which doesn't change the intersection number
by a homotopy argument), we can assume that $f$ is a {\em proper stratified Morse
function}, i.e. all critical points of $f$ with respect to $X_{\bullet}$
are {\em stratified Morse critical points} in the sense of \cite[def.5.0.2, p.277]{Sch}
and \cite[def.2.1, p.52]{GM}:
if $x\in X^{s}$ is a {\em critical} point of $f|X^{s}$, then $df_{x}$ 
is a {\em nondegenerate} covector and $f|X^{s}$ has in $x$ a classical
{\em Morse critical point} (i.e. its Hessian is nondegenerate).\\

By remark \ref{rem:musupp} and the assumption on $df(X)$ (or $-df(X)$), there 
are then only {\em finitely} many
such critical points $x$, whose ``normal Morse datum'' with respect to
${\cal F}$ are non-trivial at $df_{x}$ (or $-df_{x}$).
By \cite[lem.5.1.1, p.303]{Sch} and \cite[lem.5.1.2, p.306]{Sch} we get for $r\in [a,b[\;$
big enough
(by a {\em Mittag-Leffler argument} as in \cite[cor.6.1.2, p.430]{Sch}):
\[R\Gamma(X,{\cal F})\simeq R\Gamma(\{f\leq r\},{\cal F})\]
in the first case, and
\[R\Gamma_{c}(\{f<r\},{\cal F})\simeq R\Gamma_{c}(X,{\cal F})\]   
in the second case. By induction, it is enough to consider the case
that $f$ has at most one critical point $x\in \{f\leq r\}$, with
$x\in X^{s}$, and $f(x)<r$. Then we get by \cite[thm.5.0.1(1), p.278]{Sch} and 
\cite[lem.5.1.2, p.306]{Sch}:
\[R\Gamma(\{f\leq r\},{\cal F})\simeq NMD({\cal F},f,x)[-\tau]
\quad \text{or} \quad
R\Gamma_{c}(\{f< r\},{\cal F})\simeq NMD({\cal F},-f,x)[-\tau] \:,\]
with $\tau$ the {\em Morse index} of $\pm f|X^{s}$ (i.e. its Hessian in $x$ has
$\tau$ negative eigenvalues).
But the corresponding 
``normal Morse datum'' $NMD({\cal F},\pm f,x)$ is finitely dimensional (as explained
before). Therefore our claim follows from 
the description (\ref{eq:CC}) for $CC({\cal F})$, and the {\em local
intersection formula} (see \cite[eq.(9.5.18), p.388]{KS} and \cite[eq.(5.20), p.292]{Sch}): 
\[\sharp_{df_{x}}(\:[\;\pm \; df(X)]\cap [T^{*}_{X^{s}}X]\:) = (-1)^{\tau}\:.\]
\end{proof}

\begin{rem} We used in the above proof the fact, that we can approximate $f$
by a stratified Morse function. This follows from our assumption
about the (local) existence of a suitable {\em stratification}   
of the set of {\em degenerate} covectors. By using a proper $C^{2}$-embedding
$X\hookrightarrow \bb{R}^{N}$, we can assume $X=\bb{R}^{N}$.
Then our claim follows for example from \cite[thm.1, thm.3]{Orro}.
\end{rem}


\section{Intersection formula}

In this section we work in one of the following ``geometric categories'':
\begin{enumerate}
\item $X$ is an affine space $\bb{R}^{n}$, and $\mathcal{S}$ is 
an {\em o-minimal structure} on the real field $(\bb{R},+,\cdot)$ 
(\cite{vDr, vDrMi}, e.g. the structure of {\em semialgebraic} subsets
of real affine spaces). We could also assume, that $X$ is a real analytic
$\mathcal{S}$-manifold as in \cite[p.507-508]{vDrMi} (e.g. a
real analytic Nash manifold).
\item $X$ is a real analytic manifold, and $\mathcal{S}$ is an
{\em analytic geometric category} (\cite{vDrMi, SchVi}, e.g. the structure
of {\em subanalytic} subsets of real analytic manifolds).
\item $X$ is a real analytic Nash manifold, and $\mathcal{S}$ is a
{\em Nash geometric category} (\cite[chapter 2]{Sch} and \cite{Sch3}, e.g. the structure
of {\em locally semialgebraic} subsets of real analytic Nash manifolds).
\end{enumerate}

The subsets of $\mathcal{S}(X)$ are by definition the {\em definable} subsets 
of $X$. A {\em definable map} $f: A\to B$ between definable subsets $A\subset X$
and $B\subset X'$ (with $X,X'$ ambient manifolds as in the above cases)
is a continous map with {\em definable graph}. A complex of sheaves
${\cal F}$ on $X$ is called {\em $\mathcal{S}$-constructible}, if it is 
constructible with respect to a filtration $X_{\bullet}$:
\[\emptyset=X_{-1}
\subset X_{0}\subset \dots \subset X_{n}=X\] 
by closed {\em definable} subsets of $X$ (i.e. the cohomology sheaves of
${\cal F}$ are {\em locally constant} on all $X^{i}:=X_{i}\backslash X_{i-1}$).
By \cite[1.19, 4.8]{vDrMi} we can then assume, that ${\cal F}$ is constructible
with respect to a {\em definable Whitney b-regular $C^{p}$-stratification},
with $1\leq p < \infty$ (i.e. a definable filtration as before, with
all $X^{i}$ $C^{p}$-submanifolds of $X$, which is b-regular).

\begin{rem} For the following discussion, the reader should also compare
with \cite[sec. 10]{SchVi}. For the basic results about definable sets and 
maps, we refer to \cite{vDrMi}. Note, that their results 
about {\em analytic geometric categories} are also true
in the context of {\em Nash geometric categories} (with the obvious
modifications, as explained in 
\cite[chapter 2]{Sch}, where we developed the basic theory of
$\mathcal{S}$-constructible complexes of sheaves).
\end{rem}

Let $X_{\bullet}$ be a definable Whitney stratification. Then
\[\Lambda':= \bigcup_{i} \: \Bigl(T_{X^{i}}^{*}X \backslash \bigcup_{i\neq j} \:
cl(T_{X^{j}}^{*}X) \Bigr) \;\subset \; 
\Lambda:=\bigcup_{i} \: T_{X^{i}}^{*}X \]
are {\em definable} subsets of $T^{*}X$, and the dimension of 
the set $\Upsilon:=\Lambda\backslash \Lambda'$ of {\em degenerate} covectors
is $<d:=dim(X)$. Therefore, the set of {\em nondegenerate} covectors $\Lambda'$
is {\em dense} in $\Lambda$, and we can find (\cite[1.19, 4.8]{vDrMi})
a {\em Whitney b-regular stratification} $\Lambda_{\bullet}$ of $\Lambda$ of
the form
\[\emptyset=\Lambda_{-1}
\subset \Lambda_{0}\subset \dots \subset \Lambda_{d-1}=\Upsilon
\subset \Lambda_{d}=\Lambda \:.\]
In particular, the frontier 
$\partial \Lambda'_{j}:=cl(\Lambda'_{j}) \backslash \Lambda'_{j}$
of each connected component $\Lambda'_{j}$ of $\Lambda'$ is 
a union of strata of $\Upsilon$ (\cite[prop.4.0.2(1), p.219]{Sch}). \\

Therefore we can apply the description of the first section for 
$|CC({\cal F})|, \mu supp({\cal F})$ and $[CC({\cal F})]$,
and theorem \ref{thm:int3} implies theorem \ref{thm:int2}.\\

Now we are ready for the proof of the {\em main theorem} of this paper:
 
\begin{thm} \label{thm:int5} 
Let $f: X\to \bb{R}$ be a definable $C^{2}$-function. 
Consider on $X$ a bounded $\mathcal{S}$-constructible complex
${\cal F}$ of sheaves of vector spaces, with finite 
dimensional stalks ${\cal F}_{x}$ ($x\in X$). Suppose that the intersection
of $\sigma_{f}$ and the support 
$|CC({\cal F})|$ of the characteristic cycle of ${\cal F}$ is 
contained in a
compact definable subset $I\subset T^{*}M$, with $K:=\pi(I)\subset \{f=0\}$.
Then one has
\begin{equation} \label{eq:int1def}
\chi\bigl(R\Gamma(K,R\Gamma_{\{f\geq 0\}}{\cal F})\bigr) = 
\sharp\bigl(\; [df(X)]\cap [CC({\cal F})]\; \bigr) \:.
\end{equation}
\end{thm}

Note that $R\Gamma(K,R\Gamma_{\{f\geq 0\}}{\cal F})$ has finite dimensional
cohomology, since $K$ is a compact definable subset (compare 
\cite[chapter 2]{Sch}).

\begin{proof} 1.) Since $K$ is a compact definable subset of $X$, 
there exists by \cite[1.20, 4.22, D.19]{vDrMi} a definable $C^{2}$-function
$g: X\to \bb{R}$, with $\{g=0\}=K$. We can also assume that $g$ is 
non-negative (otherwise use $g^{2}$).\\

Choose a definable Whitney $b$-regular $C^{2}$-stratification
$X_{\bullet}$ of $X$ with the following properties (\cite[1.19, 4.8]{vDrMi}):
\begin{enumerate}
\item[a.] ${\cal F}$ is {\em constructible} with respect to $X_{\bullet}$.
\item[b.] The definable sets $\{f=0\}$ and $K$ are {\em unions of strata}.
\item[c.] The set of {\em nondegenerate} covectors (with respect to this 
stratification) is {\em dense} in all fibers of the projection 
$\pi': T^{*}_{S}M \to S$ (for all strata $S$). 
\end{enumerate}
The last property can be achieved inductively
by \cite[1.19, 4.8]{vDrMi} and the fact, that the dimension of 
the set of {\em degenerate} covectors $\Upsilon$ is $<d:=dim(X)$.
Note, that the set of nondegenerate covectors is {\em dense} in such a fiber 
$\{\pi'=x\}$, if and only if $dim(\Upsilon\cap \{\pi'=x\})<d-dim(S)$.
But the set 
\[\{x\in S\;|\; dim(\Upsilon\cap \{\pi'=x\}) \geq d-dim(S)\;\}\]
is a definable subset of $S$ with dimension $<dim(S)$ (compare
\cite[prop.1.5, p.65]{vDr}). \\

2.)  There exists a relatively compact open neighborhood $U$ of $K$ in $M$, and a $\delta_{0}>0$ 
with the following properties:
\begin{enumerate}
\item[a')] $g: \{g<\delta_{0}\}\cap U\to [0,\delta_{0}[\;$ is {\em proper}
(since $\{g=0\}=K$ is compact).
After restriction to $U$ we can and will assume $M=U$. 
\item[b')] $g$ has no {\em critical values} in $\{0<g<\delta_{0}\}$ with respect to
$X_{\bullet}$ (this follows from a') and the $C^{1}$-version of the
{\em curve selection lemma} \cite[1.17]{vDrMi}).
\item[c')] The natural morphism 
\[R\Gamma(\{g<\delta\},R\Gamma_{\{f\geq 0\}}{\cal F})\to    
R\Gamma(K,R\Gamma_{\{f\geq 0\}}{\cal F})\] 
is for {\em all} $0<\delta<\delta_{0}$
an {\em isomorphism}. This follows from the fact, that the restriction of 
$Rg_{*}({\cal F}|\{g<\delta_{0}\})$ to $[0,\delta_{0}[\;$, with $\delta_{0}$ as 
in a'), is $\mathcal{S}$-constructible 
(see \cite[thm.2.0.1, p.85]{Sch} and \cite[cor.2.2.1, p.105]{Sch}). One can also use b') and (a 
cohomologically version of) the {\em first isotopy lemma of Thom} (as in
\cite[subsec.4.1.1]{Sch}). 
\item[d')] $k\geq 0,\: 0<g(x)<\delta_{0},\: f(x)>0\:\Rightarrow \:
k\cdot dg_{x}+df_{x} \;\not\in \Lambda$.
\end{enumerate}
One shows the last property indirectly (compare \cite[p.386-387]{KS}), and
note that
\[\{(x,k)\;|\;k\geq 0,\: 0<g(x)<\delta_{0},\: f(x)>0,\:
k\cdot dg_{x}+df_{x} \;\in \Lambda\;\}\]
is a {\em definable} subset of $X\times \bb{R}$:

Otherwise there exists (\cite[1.17]{vDrMi}) a stratum $X^{s}$ and a 
$C^{1}$-curve $\gamma: [0,\delta_{0}[\;\to X\times \bb{R}^{2}$, with
$\gamma(t)=:(x(t),\alpha(t),\beta(t))$ such that $\alpha(t)\geq 0,\;
\beta(t)>0,\; g(x(0))=0,\; f(x(t))>0$ for $t>0$, and 
\[\alpha(t)\cdot dg_{x(t)} + \beta(t)\cdot df_{x(t)}\;\in T^{*}_{X^{s}}X
\quad \text{for $t>0$.} \] 
But this implies $f(x(0))=0$, since $\{g=0\}=K\subset \{f=0\}$, and
\[\alpha(t)\cdot d/dt\;\bigl(g(x(t))\bigr) + 
\beta(t)\cdot d/dt\;\bigl(f(x(t))\bigr)\:\equiv \;0\:,\]
since $T^{*}_{X^{s}}X$ is a {\em conic isotropic}
submanifold of $T^{*}X$ (cf. \cite[p.483]{KS}).

But by the {\em monotonicity theorem} \cite[4.1]{vDrMi} we have that
$f\circ x$ and $g\circ x$ are {\em constant or strictly monotonic} on 
$\;]0,a[\;$ for $a>0$ small enough. By the assumptions, this implies
the contradiction $d/dt\;\bigl(g(x(t))\bigr)\geq 0$ and
$d/dt\;\bigl(f(x(t))\bigr)>0$ for $0<t<a$. \\

3.) Fix a $\delta$ with $0<\delta<\delta_{0}$, and consider the
inclusion $j: \{g<\delta\}\to X$. By the choice of $X_{\bullet}$ and b')
we can apply theorem \ref{thm:int3}. Let ${\cal F}':=Rj_{*}j^{*}{\cal F}$,
which is {\em constructible} with respect to the induced definable b-regular 
stratification of $\{g\leq \delta\}$. By the {\em curve selection lemma}
we can find $\epsilon>0$ so that $f$ has in $\{0<|f|\leq\epsilon\}$ {\em no 
critical points} with respect to this induced stratification.\\

Note that $supp({\cal F}')\subset
\{g\leq \delta\}$ is {\em compact}. Therefore we get by \cite[lem.5.0.1, p.276]{Sch}
and \cite[lem.5.1.1(2), p.303]{Sch} (compare also with \cite[cor. 5.4.19.(ii), p.239]{KS}):
\begin{equation} \label{eq:compact}
R\Gamma(X,R\Gamma_{\{f\geq 0\}}{\cal F}') \simeq  
R\Gamma_{c}(\{f>-\epsilon\},{\cal F}')\:.
\end{equation}
Together with c') and the general isomorphism
\[R\Gamma(X,R\Gamma_{\{f\geq 0\}}{\cal F}') \simeq  
R\Gamma(\{g<\delta\},R\Gamma_{\{f\geq 0\}}{\cal F})\:,\]
this implies the isomorphism
\begin{equation} \label{eq:compact2}
R\Gamma(K,R\Gamma_{\{f\geq 0\}}{\cal F}) \simeq  
R\Gamma_{c}(\{f>-\epsilon\},{\cal F}')\:.
\end{equation}

Choose an $a\in \bb{R}$ with $supp({\cal F}')\subset \{f<a\}$,
and consider $h:=-f: U:=\{-\epsilon<f<a\}\to [-a,\epsilon[\;$.
Note that $h|supp({\cal F}'|U)$ is {\em proper}.
Since $f$ has in $\{0<|f|\leq\epsilon\}$ no 
critical points with respect to the induced stratification
of $\{g\leq\delta\}$, we get
\begin{equation} \label{eq:compact3}
df(U)\cap \mu supp({\cal F}'|U) \subset T^{*}U|_{\{f\geq 0\}}\:.
\end{equation}
In particular, $df(U)\cap \mu supp({\cal F}'|U)\subset
df(supp({\cal F}')\cap \{f\geq 0\})$ is a {\em compact} subset of 
$T^{*}U$. Therefore we can apply theorem \ref{thm:int4}.2 
to the manifold $U$, the function $h$ and the sheaf complex ${\cal F}'|U$,
and get the intersection formula:
\begin{equation} \label{eq:compact4}
\chi(R\Gamma_{c}(\{f>-\epsilon\},{\cal F}')) = 
\sharp(\:[df(U)]\cap [CC({\cal F}'|U)]\:)\:.
\end{equation}

By equation (\ref{eq:compact2}), the proof of the theorem is complete,
if we show
\[|CC({\cal F}')| \cap \{\;df_{x}\;|\; g(x)=\delta,\;f(x)>-\epsilon\;\}    
= \emptyset\:.\]
For $-\epsilon<f(x)<0$, we get $df_{x}\not\in |CC({\cal F}')|$ by 
(\ref{eq:compact3}).\\

Assume $df_{x}\in |CC({\cal F}')|$, with $f(x)\geq 0$ and $g(x)=\delta$.
By theorem \ref{thm:int3} we get
$df_{x} = \omega -c\cdot dg_{x}$, with $c\geq 0,\; \pi(\omega)=x$ and
$\omega\in |CC({\cal F})|$. In the case $f(x)>0$ we get a contradiction
to d'): 
\[c\cdot dg_{x} + df_{x} = \omega \in |CC({\cal F})|\subset \Lambda \;.\]
So we can assume $f(x)=0$. In the case $c=0$ we get 
\[df_{x}= \omega \in |CC({\cal F})|\;.\]
But this is by $g(x)=\delta >0$ impossible, since
the intersection of $\sigma_{f}$ and $|CC({\cal F})|$ is 
contained in the compact subset $I\subset T^{*}M$, 
with $\{g=0\}=K:=\pi(I)$. In the case $c>0$ we get 
\[dg_{x}=(\omega -df_{x})/c \in \Lambda\;,\]
since $|CC({\cal F})|\subset \Lambda$, 
and $df_{x}\in \Lambda$ for $f(x)=0$ (because $\{f=0\}$ is by b. a union of
strata of $X_{\bullet}$). But this is a contradiction to b').
\end{proof}


\end{document}